\documentclass[reqno]{article}
\usepackage{amsfonts}
\usepackage{amsmath}
\usepackage{amsthm}
\usepackage{amssymb}
\usepackage{oldgerm}
\usepackage{fancyhdr}
\usepackage{fancybox}
\usepackage{mathrsfs}
\usepackage{epsfig}
\def\phi{\varphi }
\swapnumbers

\theoremstyle{plain}
\newtheorem{theorem}{Theorem}[section]
\newtheorem{corollary}[theorem]{Corollary}
\newtheorem{lemma}[theorem]{Lemma}
\newtheorem{proposition}[theorem]{Proposition}
\theoremstyle{definition}
\newtheorem{definition}[theorem]{Definition}
\theoremstyle{remark}

\newtheorem{example}[theorem]{Example}

\numberwithin{equation}{section}

\begin{document}
\title{Positive convolution structure for a class of
Heckman-Opdam hypergeometric functions of type BC}
\author{Margit R\"osler\\
Institut f\"ur Mathematik, TU Clausthal\\
Erzstr. 1\\
D-38678 Clausthal-Zellerfeld, Germany\\
roesler@math.tu-clausthal.de}
\date{}
\maketitle

\begin{abstract}
In this paper, we derive explicit product formulas and positive convolution structures 
for three continuous classes of Heckman-Opdam hypergeometric functions of type $BC$. For specific discrete series of  multiplicities
these hypergeometric functions occur as the spherical functions of non-compact Grassmann manifolds $G/K$ over one of the skew fields $\mathbb F= \mathbb R, \mathbb C, \mathbb H.$ 
We write the product formula of these spherical functions  in an explicit form which allows analytic continuation with respect to the parameters. In each of the three cases, we obtain a  series of hypergroup algebras
 which include the commutative convolution algebras of $K$-biinvariant functions 
on $G$ as special cases. The characters are given by the associated hypergeometric functions.

\end{abstract}




\section{Introduction}
There is a well-established theory of hypergeometric functions associated with root systems due to Heckman, Opdam and Cherednik which generalizes and completes the theory of 
spherical functions on Riemannian symmetric spaces in many respects; see \cite{O1}, \cite{HS}, \cite{O2}, \cite{S} as well as the literature cited there. In rank one, i.e. for root systems of type $BC_1$,  these hypergeometric functions are known as  Jacobi functions and were studied by Flensted-Jensen and Koornwinder in a series of papers  in the 1970ies. A comprehensive exposition is given in \cite{K}. In generalization of the one-variable case, hypergeometric functions associated with root systems are indexed by continuous parameters (the multiplicities) on a given root system. They build up the solutions of the joint eigenvalue problem for an associated system of commuting differential operators which
generalize the radial parts of all invariant differential operators on a Riemannian symmetric space $G/K$  of the non-compact type. In such geometric cases, the root system and multiplicity function are given in terms of the root space data of $(G, K)$. 
In fact, the harmonic analysis associated with such hypergeometric functions  is only the Weyl-group invariant part of a more general harmonic analysis
associated with a commuting family of differential-reflection operators of Dunkl type, the so-called Cherednik operators. The associated integral transform, which generalizes
the spherical transform on symmetric spaces, is studied in detail in \cite{O1}. There are, in particular, a Paley-Wiener theorem and a Plancherel theorem established for this transform. In the geometric cases
$(G,K)$ is a Gelfand pair, and the corresponding spherical functions satisfy a product formula
which is intimately connected to the harmonic analysis on the commutative algebra of $K$-biinvariant 
measures on $G$. In the rank one case, a positive product formula and harmonic analysis for Jacobi functions associated with general non-negative multiplicities was established by Flensted-Jensen and Koornwinder, see \cite{K}.
However, apart from theses cases, the existence of a positive product formula for multivariable hypergeometric functions and a positivity-preserving convolution which would allow for 
a general $L^p$-theory are still open in general. 

A natural idea to extend the convolution from particular geometric cases
to general multiplicities is analytic continuation of the product formula with respect to the multiplicities. There are only three classes of geometric cases with an infinite discrete series of multiplicities when the rank is fixed, namely the 
non-compact Grassmann manifolds $SO_0(p,q)/SO(p)\times SO(q), SU(p,q)/S(U(p)\times U(q))$ and  $Sp(p,q)/Sp(p)\times Sp(q).$ 
Their real  rank is $q$ and the spherical functions are hypergeometric functions of type $BC$ with multiplicities depending on $p$. In the present paper, we carry out the interpolation program in these cases. We 
give an explicit product formula for the spherical functions which allows analytic extension with respect to the multiplicity parameter $p$. This yields a product formula for three continuous classes of hypergeometric functions of type $BC$ interpolating the group cases. 
Based on the product formula, we obtain a complete picture 
of harmonic analysis  within the framework of commutative hypergroups on the associated Weyl chamber. In particular, the hypergeometric transform becomes an interpretation as a hypergroup Fourier transform. 

The paper is organized as follows: In Section \ref{spherical}, we calculate the product formula for the spherical  functions on the Grassmann manifolds. Section \ref{HO_theory} gives a short account on Heckman-Opdam theory as well as the identification of the spherical functions on Grassmann manifolds as hypergeometric functions of type $BC_q$. The extension of the product formula to a continuous range 
of multiplicities interpolating the dimension parameter $p$ is carried out in Section \ref{prod_formula}, and Section \ref{hypergroups} is devoted to the study of the associated hypergroup algebras on the 
Weyl chamber. A central part of this section is the characterization of the bounded multiplicative functions which generalizes well-known results for spherical functions. The reasoning here is, however, not based on an integral representation but on exponential bounds for the Heckman-Opdam hypergeometric functions and  their generalized Harish-Chandra expansion.

\section{Spherical functions on Grassmann manifolds and their product formula}\label{spherical}

We consider the Grassmann manifolds $G/K$ where $G$ is one of the  indefinite orthogonal, unitary or symplectic  groups
$ SO_0(p,q),\, SU(p,q)$  or $Sp(p,q)$ with maximal compact subgroup  $K= 
SO(p)\times SO(q), \, S(U(p)\times U(q))$ or $Sp(p)\times Sp(q),$ respectively. For a unified point of view we also consider  $K$ as subgroup of $U(p;\mathbb F)\times U(q,\mathbb F),$ where   $U(p;\mathbb F)$ is the unitary group
over  $\mathbb F = \mathbb R,  \mathbb C$ or  $\mathbb H$. 
In the same way  $G$ is a subgroup of  the indefinite unitary group $U(p,q;\mathbb F),$ which is the isometry group for the quadratic form
\[ |x_1|^2 + \ldots +|x_p|^2 - |x_{p+1}|^2 -\ldots - |x_{p+q}|^2\]
on $ \mathbb F^{p+q}$. To avoid exceptions which will be irrelevant lateron, we shall exclude the case $p=q$ and assume that $p>q\geq 1$.

It is well known that $(G,K)$ is a Gelfand pair (this follows from Corollary 1.5.4. of \cite{GV}). The spherical functions of this pair are characterized as the non-zero $K$-biinvariant
continuous functions $\phi:G\to \mathbb C$ which satisfy the product formula 
\begin{equation}\label{prodformula}
\phi(g)\phi(h) = \int_K \phi(gkh)dk \quad \text{ for all }\, g,h \in G
\end{equation}
where $dk$ denotes the normalized Haar measure of $K$. 
This means that the space of continuous, $K$-biinvariant compactly supported functions on $G$ 
is a commutative subalgebra of the convolution algebra $C_c(G)$. 
The space $C_c(G//K)$ on the double coset space $G//K$ therefore inherits the structure of a commutative topological algebra. 
 The spherical functions of $(G,K)$ provide exactly the non-zero continuous characters of this algebra, via $\, 
f\mapsto \int_G f(x)\phi(x)dx.$ 

\medskip
To make the product formula explicit, we 
recall the $K\!AK$-decomposition of $G$.  Let $\frak g$ and $\frak k$ denote the Lie algebras of $G$ and $K$. $\frak  g$  has the Cartan decomposition $\frak g = \frak k \oplus \frak p$ with $\frak p$ consisting of the $p+q$-block matrices
\[ \begin{pmatrix} 0& X\\
    X^* & 0\end{pmatrix}, \quad X\in M_{p,q}(\mathbb F).\]
Let 
$\frak a$ be a maximal abelian subalgebra of $\frak p.$ Then $G=K\!AK$ with  $A=\exp \frak a.$
The spherical functions of $(G,K)$ are therefore determined by their values on $A$. 
Actually, they are already determined by
their values on the topological closure $\overline{A_+} = \exp(\overline {\frak a_+})$  if $\frak a_+$  is the positive Weyl chamber associated with an (arbitrary) choice of positive roots
within the restricted root system $\Delta= \Delta(\frak a, \frak g)$ of $\frak g$ with respect to $\frak a$. 
We may choose for  $\mathfrak a$  the set of all matrices
$H_t\in M_{p+q}(\mathbb F)$ of the form 
\[ H_t = \begin{pmatrix} 0_{p\times p}& \begin{matrix}\underline t \\
                            0_{(p-q)\times q}\end{matrix}\\
\begin{matrix} \,\underline t& 0_{q\times(p-q)}\end{matrix}& 0_{q\times q}\end{pmatrix} \]
where $\, \underline t := \text{diag}(t_1, \ldots, t_q)$ is the $q\times q$ diagonal matrix corresponding to $t=(t_1,  \ldots, t_q)\in \mathbb R^q$ (here $\mathbb R$ is considered as a subfield of $\mathbb C$ and $\mathbb H$ in the usual way).
The real rank of $G$ is $q$, and the restricted root system $\Delta= \Delta(\frak a, \frak g)$ is 
of type $BC_q$ with the understanding  that zero is allowed as a multiplicity on the long roots.   In this way the limiting case $B_q$, which occurs 
for $\mathbb F=\mathbb R$, is included.
We identify $\frak a$ with $\mathbb R^q$ via $H_t \mapsto t$, where the coordinates are with respect to the standard basis ${e_1, \ldots, e_q}$ of $\mathbb R^q$. Then the Killing form on $\frak a$ becomes the standard Euclidean inner product on $\mathbb R^q$. Here is a comprehensive table of the roots $\alpha$ and their (geometric) multiplicities $m(\alpha)$, that is the dimensions of the corresponding root spaces; c.f. Table 9 of \cite{OV}. The constant $d$ denotes the dimension of $\mathbb F$ as an $\mathbb R$-vectorspace, i.e. $d=1,2,4$ for $\mathbb F= \mathbb R,\mathbb C, \mathbb H$. 

\begin{equation}\label{tabelle}
\begin{tabular}{|c|c|}\hline
root $\alpha$ &  multiplicity  $m(\alpha)= m_{p,d}(\alpha)$  \\ \hline
$\alpha(t) = \pm t_i\,;\, \,1\leq i \leq q$ & $d(p-q)$ \\ \hline
$\alpha(t) = \pm 2 t_i\,; \, 1\leq i \leq q$ & $d-1$ \\ \hline
$\alpha(t) = \pm t_i\pm t_j\,;\, 1\leq i < j \leq q $ & $d$ \\ \hline
\end{tabular}
\end{equation}

\medskip
Thanks to our restriction $p>q$, the Weyl group of $(\frak a,\frak g)$ is the hyperoctahedral group in all cases, and as a Weyl chamber we may choose 
\[\frak a_+ := \{H_t:\, t=(t_1, \ldots t_q)\in \mathbb R \text{ with }\,  t_1 >t_2 > \ldots > t_q >0 \}.\]
In our identification of $\frak a$ with $\mathbb R^q$,  the closed chamber $\overline{\frak a_+}$ corresponds to the set
\[ C := \,\{ t\in \mathbb R^q: \, t_1 \geq t_2 \ldots \geq  t_q \geq 0\}.\]
A short calculation gives
\[\overline{A_+} = \left\{ a_t = \, \begin{pmatrix} \cosh \underline t & 0_{q\times (p-q)} & \sinh \underline t\\
0_{(p-q)\times q}& I_{p-q}& 0_{(p-q)\times q} \\
\sinh \underline t & 0_{q\times (p-q)} & \cosh \underline t
            \end{pmatrix} \in M_{p+q}(\mathbb F)\,:\,\, t\in C\right\}.\]

\smallskip\noindent
Consider now
\[ g = \begin{pmatrix} u & 0 \\
0 & v \end{pmatrix} a_t \begin{pmatrix} \widetilde u & 0\\ 0 & \widetilde v
\end{pmatrix} \in K\!a_t K.\]
To obtain $t$ back from $g$, we write $g$ in $p\times q$ block notation as 
\[ g = \begin{pmatrix} A(g) & B(g)\\
C(g) & D(g)
       \end{pmatrix}.\]
A short calculation gives 
\begin{equation}\label{dg}
 D(g) = v \cosh\underline t \, \widetilde v. 
\end{equation}
Let $\text{spec}_s(x)$ denote the singular spectrum of $x\in M_q(\mathbb F),$ that is 
\[\text{spec}_s(x)= \sqrt{\text{spec}(x^*x)} = (\lambda_1,\ldots, \lambda_q)\in \mathbb R^q\]
with the singular values $\lambda_i$ of $x$ ordered by  size: $\lambda_1 \geq \ldots \geq \lambda_q\geq 0.$ 
Equation \eqref{dg} shows that the singular spectrum of $D(g)$ is given by
$\,  \text{spec}_s (D(g)) \, = (\cosh t_1, \ldots \cosh t_q) = : \cosh t \,$. Therefore
\begin{equation}\label{parameterident}
 t = \text{arcosh}\bigl(\text{spec}_s(D(g)\bigr) \quad \text{ for each }\, g\in K\!a_t K, \, t\in C
\end{equation}
where arcosh is also taken componentwise.
(Observe that $\,D(g) \geq I_q$ and therefore all its singular values are $\geq 1$).

Let us now evaluate the product formula \eqref{prodformula} for the spherical functions of $(G,K)$ explicitly. 
As spherical functions are $K$-biinvariant, it suffices
to calculate the product formula for arguments $\,g= a_t, h=a_s \in \overline{A_+}$. Write $a_t\in  \overline{A_+}$ in $p\times q$-block notation:
\[a_t = \begin{pmatrix} A_t & B_t\\
   C_t & D_t
      \end{pmatrix}.\]
Then for $a_t, a_s\in \overline{A_+}$ and $k= \begin {pmatrix} u & 0\\ 0 & v\end{pmatrix} \in K\,$ we obtain
\[ a_tka_s = \begin{pmatrix} * & *\\
* & C_tuB_s +D_tvD_s
         \end{pmatrix}
\]
and therefore 
\[ D(a_tka_s) = C_tuB_s +D_tvD_s \,=\,\bigl(\sinh\underline t \,\vert\, 0\bigr) u \begin{pmatrix} \sinh\underline s \\
 0 \end{pmatrix} \,+\, \cosh\underline t \,v \cosh\underline s.\]
With the block matrix 
\[ \sigma_0 := \begin{pmatrix} I_q\\ 0\end{pmatrix} \in M_{p,q}(\mathbb F)\]
this can be written as 
\[ D(a_tka_s) = \,\sinh\underline t \,\sigma_0^* u \sigma_0\sinh\underline s \,+\,\cosh\underline t \,v \cosh\underline s.\] 
Notice that $\,\sigma_0^* u \sigma_0\in M_q(\mathbb F)$ is a truncation of $u$ given by the upper left $q\times q$-block of $\sigma$. 

\smallskip\noindent
Let  $\phi$ be a spherical function of $(G,K)$ and put $\, \widetilde\phi(t) := \phi(a_t)\,$ for $t\in C$. 
Then according to formula \eqref{parameterident} it satisfies
\begin{equation}\label{prodformel1}
 \widetilde\phi(t)\widetilde\phi(s) \, = \, \int_K \widetilde\phi\bigl(\text{arcosh}\bigl(\text{spec}_sD(a_t ka_s)\bigr)\bigr) dk.\end{equation}
In order to achieve a simplification of this formula  we first extend the integral over $K$ to an integral over
 $U(p;\mathbb F)\times U_0(q;\mathbb F) =: K_0$, where $U_0(q;\mathbb F)$ denotes the connected component of the identity in $U(q;\mathbb F)$.
If $\mathbb F=\mathbb H$ then $K=K_0$, 
but in the other cases  $K$ is a proper normal subgroup of $K_0.$ More precisely, let $\mathbb T:=\{z\in \mathbb F: |z|=1\}$ and $H$ the group of diagonal matrices 
$H=\{ d_z: z\in \mathbb T\} \subset M_{p+q}(\mathbb F)$ where the diagonal entries of $d_z$ are equal $1$ apart from the entry in position $(p,p)$, which is $z$. Then $K_0 = H\ltimes K\cong \mathbb T\ltimes K.$
Suppose $f$ is a continuous function on $K_0$ of 
the form 
\[ f(k_0)= \widetilde f(\sigma_0^*u\sigma_0,v) \quad \text{for } \,k_0=\begin{pmatrix}u & 0\\
                                                                                                                   0 & v
                                                        \end{pmatrix}.\]
Then $\,f(d_zk) = f(k) \,$ for all $z\in \mathbb T$ and $k\in K$ and thus 
by Weyl's formula, 
\[ \int_{K_0} f(k_0)dk_0\,=\, \int_{\mathbb T} \bigl(\int_{K} f(d_zk)dk\bigr)dz\,= \, 
\int_K f(k)dk\]
where on each of the involved groups, integration is with respect to the normalized Haar measure.
Thus
\[\widetilde\phi(t)\widetilde\phi(s)\,=\, 
 \int_{U(p,\mathbb F)}\int_{U_0(q,\mathbb F)} \widetilde\phi\bigl(\text{arcosh}(\text{spec}_s(\sinh \underline t \,\sigma_0^* u \sigma_0\sinh\underline s \,+\,\cosh\underline t \,v \cosh\underline s))\bigr) dudv\]
with $du$ and $dv$ the normalized Haar measures on  $U(p,\mathbb F)$ and $U_0(q,\mathbb F) $ respectively.
Here the integrand depends only on $v$ and the truncation $\sigma_0^*u\sigma_0$, which  is contained in the closure of the ball
\[ B_q := \{w\in M_q(\mathbb F): w^*w < I\}.\]

Under the assumption $p\geq 2q$ this situation is covered by the following reduction lemma, which is a consequence of Corollary 3.3. of \cite{R}. Let
\[ \gamma := d(q-\frac{1}{2}) +1\]
and for $\mu\in \mathbb C$ with $\text{Re}\, \mu >\gamma -1,$ put 
 \begin{equation}\label{kappa} \kappa_\mu = \int_{B_q} \Delta(I-w^*w)^{\mu-\gamma} dw.
 \end{equation}

Here $\Delta(x)$ denotes the determinant of $x\in M_q(\mathbb F)$, which is 
defined as the usual determinant for $\mathbb F=\mathbb R$ or $\mathbb C$, while for $\mathbb F=\mathbb H$ we choose the Dieudonn\'e determinant, i.e. $\Delta(x) = (\det_{\mathbb C}(x))^{1/2}$ when $x$ is considered as a complex matrix in the usual way. 

\begin{lemma} Suppose that $p\geq 2q$. 
Then for continuous $f: \overline B_q\to \mathbb C$, 
\[ \int_{U(p,\mathbb F)}f(\sigma_0^*u\sigma_0) du\,=\, \frac{1}{\kappa_{pd/2}} \int_{B_q}f(w) \Delta(I-w^*w)^{pd/2-\gamma}dw.\]
\end{lemma}

\begin{proof} Consider the action of  the unitary group $U(p,\mathbb F)$ on $M_{p,q}(\mathbb F)$ by left multiplication, $\,(u,x) \mapsto\, ux.\,$
The orbit of the matrix $\sigma_0$ under this action is the Stiefel manifold
\[ \Sigma_{p,q} = \{ x\in M_{p,q}(\mathbb F): x^*x = I_q\}.\]
Consider further  the map $U(p,\mathbb F) \to \Sigma_{p,q}, \, u\mapsto u\sigma_0$. The image measure of $du$ under this map 
coincides with the normalized $U(p,\mathbb F)$-invariant measure $d\sigma$ on $\Sigma_{p,q}$.
Therefore
\[  \int_{U(p,\mathbb F)}  f(\sigma_0^*u\sigma_0) du \, = \int_{\Sigma_{p,q}}  f(\sigma_0^*\sigma) d\sigma.\]
But $\sigma_0^*\sigma$  is the $q\times q$ matrix given by the first $q$ rows of $\sigma$ only.
According to Corollary 3.3. of \cite{R}, 
\begin{equation}\label{splitting}\int_{\Sigma_{p,q}}  f(\sigma_0^*\sigma) d\sigma = \frac{1}{\kappa_{pd/2}} \int_{B_q}f(w) \Delta(I-w^*w)^{pd/2-\gamma}dw,\end{equation}
which finishes the proof.
\end{proof}

\noindent
We thus obtain 

\begin{proposition}\label{prodspher}
Suppose that $p\geq 2q$. Then the spherical functions $\,\widetilde\varphi(t) = \varphi(a_t)$ satisfy the 
product formula 
\begin{align*}\widetilde\phi(t)\widetilde\phi(s) \,
= \frac{1}{\kappa_{pd/2}} \int_{B_q}\int_{U_0(q,\mathbb F)} \widetilde\phi\bigl({\rm arcosh}({\rm spec}_s(\sinh\underline t\, w& \sinh \underline s +\, \cosh\underline t \,v \cosh\underline s))\bigr)\cdot\\
& \cdot \Delta(I-w^*w)^{pd/2-\gamma} dv dw.
\end{align*}
\end{proposition}

Notice that the dependence  on $p$ now occurs only in the density, not in the domain of integration.

\section{The spherical functions as $BC_q$-hypergeometric functions}\label{HO_theory}

In this section, we first provide the necessary background on hypergeometric functions associated with root systems.  For an introduction to the subject, we refer to \cite{O1}, \cite{O2}  and part I of \cite{HS}. In a second part, we
identify the spherical functions on Grassmann manifolds within this framework.  

 Let $\frak a$ be a finite-dimensional Euclidean space with inner product $\langle\,.\,,.\,\rangle$ which is extended to a complex bilinear form on the complexification $\frak a_\mathbb C$ of $\frak a$. We identify $\frak a$ with its dual space $\frak a^* = \text{Hom}(\frak a, \mathbb R)$ via the 
given inner product. Let $R\subset \frak a$ be a (not necessarily reduced) root system and let  $W$ be the Weyl group of $R$.
 For $\alpha\in R$ we write $\alpha^\vee = 2\alpha/\langle\alpha,\alpha\rangle$ and denote by $\sigma_\alpha(x) = x - \langle x,\alpha^\vee\rangle \alpha\,$ the orthogonal reflection in the hyperplane perpendicular to $\alpha$. 

A multiplicity function on $R$ is a function $k: R\to \mathbb C$ which is $W$-invariant, i.e. $k(w\alpha) = k(\alpha)$ for all $\alpha\in R$. We denote by $K$ the vector space of multiplicity functions on $R$ and fix
a positive subsystem $R_+$ of $R$.    For $k\in K$ we put 
\[\rho(k) := \frac{1}{2}\sum_{\alpha\in R_+} k(\alpha)\alpha.\]
The Cherednik operator in direction $\xi\in \frak a$ is the differential-reflection operator on $\frak a_{\mathbb C}$ defined by
\[ T_\xi(k)  =\partial_\xi + \sum_{\alpha\in R_+} k(\alpha) \langle\alpha,\xi\rangle \frac{1}{1-e^{-\alpha}}(1-\sigma_\alpha) - \langle\rho(k),\xi\rangle\]
where $\partial_\xi$ is the usual directional derivative  and  $e^\lambda(\xi) := e^{\langle\lambda,\xi\rangle}$
for $\lambda,\xi \in \frak a_{\mathbb C}$. For fixed multiplicity $k$, the operators $\{T_\xi(k),\, \xi\in \frak a_\mathbb C\}$ commute.
Therefore the assignment $\xi\mapsto T_\xi(k)$ uniquely extends to a homomorphism on the symmetric algebra $S(\frak a_{\mathbb C})$ over $\frak a_{\mathbb C}$, which my be identified with the algebra of complex polynomials on $\frak a_{\mathbb C}$. The differential-reflection operator which in this way corresponds to $p\in S(\frak a_{\mathbb C})$ 
will be denoted by $T(p,k)$.
Let $S(\frak a_{\mathbb C})^W$ denote the subalgebra of $W$-invariant elements in $S(\frak a_{\mathbb C})$.  Then for each  $p\in S(\frak a_{\mathbb C})^W$, the Cherednik operator $T(p,k)$ coincides
with a  $W$-invariant differential operator on $C^\infty(\frak a)^W$, the $W$-invariant functions from $C^\infty(\frak a)$.
The following theorem establishes hypergeometric functions associated with root systems. It was proved  by Heckman and Opdam in a series of papers, see \cite{HS} as well as \cite{O1}.

\begin{theorem}
There exists an open regular set $K^{reg}\subseteq K$ with
$\{k\in K:\text{Re}\, k \geq 0\} \subseteq K^{reg}$  such that for each $k\in K^{reg}$ and each spectral parameter $\lambda\in \frak a_{\mathbb C}$, the hypergeometric system
\begin{equation}\label{hypersystem}
T(p,k) f = p(\lambda) f \quad \forall p \in S(\frak a_{\mathbb C})^W
\end{equation}
has a unique $W$-invariant solution $f(t) = F(\lambda,k;t)$ which is analytic on $\frak a$ and satisfies
$f(0)=1.$ 
Moreover, there is a $W$-invariant tubular neighborhood $U$ 
of $\frak a$ in $\frak a_{\mathbb C}$ such that $F$ extends  to a (single-valued) holomorphic function on $ \frak a_{\mathbb C}\times K^{reg}\times U$, which is called the hypergeometric function associated with $R$. 
$F(\lambda,k;t)$ is $W$-invariant both in $\lambda$ and $t$. 
\end{theorem}

Suppose that $k$ is real. Then for $W$-invariant polynomials $p$ with real coefficients, we have
\[ T(p,k) \overline{F(\lambda,k;\,.\,)} \, =\, p(\overline \lambda) \overline{F(\lambda,k;\,.\,)}\]
which shows that 
\begin{equation}\label{realF}
\overline{F(\lambda,k;t)} =\, F(\overline\lambda,k;t) \quad\forall\, t\in \frak a.
\end{equation}
The uniqueness of the solution to the hypergeometric system also implies the equivalence
\[F(\lambda,k;\,.\,) = F(\lambda^\prime,k;\,.\,)\, \Longleftrightarrow\, \lambda^\prime\in W.\lambda\]

Let $C_c^\infty(\frak a)^W$ denote the  $W$-invariant functions from $C_c^\infty(\frak a)$. 
The hypergeometric transform of $f\in C_c^\infty(\frak a)^W$ is defined by
\[ \mathcal Ff(\lambda) = \int_{\frak a} f(t) F(-\lambda, k; t) d\omega(t)\]
where the measure $\omega = \omega_k$ on $\frak a$ is given by 
\begin{equation}\label{weight} d\omega(t) = \prod_{\alpha\in R} \vert e^{\langle\alpha,t\rangle/2} - e^{- \langle\alpha,t\rangle/2}                       \vert^{k(\alpha)}dt\end{equation}
($dt$ denotes the Lebesgue measure on $\frak a$). There are Paley-Wiener and Plancherel theorems for this transform which are obtained  by Weyl-group symmetrization of the (non-symmetric) Cherednik transform studied in  \cite {O1}; see also \cite{O2}.
Define the measure  $\nu=\nu_k$ on $i\frak a$ by
\[d\nu(\lambda) = \frac{1}{|c(\lambda,k)|^2}d\lambda\]
where $d\lambda$ denotes the Lebesgue measure on  $i\frak a $  and $c(\,.\,,k)$ is the $c$-function on $\frak a_{\mathbb C}$,
\begin{equation}\label{c_function} c(\lambda,k) = \prod_{\alpha\in R_+} \frac{\Gamma(\langle\lambda,\alpha^\vee\rangle + \frac{1}{2}k(\frac{\alpha}{2}))}{\Gamma (\langle\lambda,\alpha^\vee\rangle + \frac{1}{2}k(\frac{\alpha}{2}) + k(\alpha))}\cdot \prod_{\alpha\in R_+} \frac{\Gamma (\langle\rho(k),\alpha^\vee\rangle + \frac{1}{2}k(\frac{\alpha}{2}) + k(\alpha))}{\Gamma(\langle\rho(k),\alpha^\vee\rangle + \frac{1}{2}k(\frac{\alpha}{2}))}\end{equation}
with the convention that $k(\frac{\alpha}{2}) = 0$ if $\frac{\alpha}{2}\notin R$.

\begin{theorem}\label{Plancherel} (\cite{O1}, Theorems 8.6 and 9.13)
\begin{enumerate} 
\item[\rm{(1)}] The hypergeometric transform $\mathcal F$ is an isomorphism from $C_c^\infty(\frak a)^W$ onto the $W$-invariant Paley-Wiener space $PW(\frak a_\mathbb C)^W$, where $ PW(\frak a_\mathbb C)$ consists of all holomorphic functions $f$ on $\frak a_\mathbb C$ satisfying the growth condition
\[ \exists\, R>0, \,\, \forall N\in \mathbb N: \,\, \sup_{\lambda\in \frak a_{\mathbb C}}(1+|\lambda|)^N e^{-R|\text{Re}\, \lambda|}|f(\lambda)| \,<\, \infty.\]
The inverse of $\,\mathcal F: C_c^\infty(\frak a)^W \,\to \, PW(\frak a_\mathbb C)^W$ is given by
\[\mathcal F^{-1}h(t) = \int_{i\frak a} h(\lambda) F(\lambda,k;t) d\nu(\lambda).\]
\item[\rm{(2)}] 
Let $f, g \in C_c^\infty(\frak a)^W$ and let $\,\frak a_+\,$ be the  Weyl chamber of $W$ corresponding to $R_+$. Then
\[ \int_{\frak a_+} f(t)\overline{g(t)} d\omega(t) \,=\,c\int_{i\frak a_+} \mathcal F f(\lambda) \overline{\mathcal  F g(\lambda)} d\nu(\lambda)\]
where  $c>0$ is a normalization constant. 
\end{enumerate}
 \end{theorem}

\noindent
According to Propos. 6.1 of \cite{O1},
\[ |F(\lambda,k;t)| \leq |W|^{1/2}\cdot e^{|\text{Re}\lambda||t|} \quad \text{ for }\, t\in \frak a, \lambda \in \frak a_{\mathbb C}. \]
Thus for $f\in C_c^\infty(\frak a)^W$ and fixed $s\in \frak a$, the function $\, \lambda \mapsto \mathcal Ff(\lambda)F(\lambda,k;s)$ belongs to $PW(\frak a_\mathbb C)^W$, and we obtain the following

\begin{corollary}\label{translate}
For $s\in \frak a$ and $f\in C_c^\infty(\frak a)^W,$ the generalized translate 
\[ \tau_sf(t) := \int_{i\frak a} \mathcal Ff(\lambda) F(\lambda,k;s)F(\lambda,k;t) d\nu(\lambda)\]
again belongs to $C_c^\infty(\frak a)^W.$ Moreover, 
\[ \mathcal F(\tau_s f)(\lambda) =\, F(\lambda,k;s)\mathcal Ff(\lambda). \]
\end{corollary}

Let us now turn to the spherical functions on the Grassmann manifolds  $G/K$. 
They are  identified with hypergeometric functions of type $BC_q$, as follows: Consider
$\frak a= \mathbb R^q$ with the standard inner product $\langle \,.\,,.\,\rangle$ and regard
the restricted root system of  $G/K$ as a subset of $\mathbb R^q$ as described in Section  \ref{spherical}.
With our convention including the case $\mathbb F=\mathbb R$,  it is given by
\[ BC_q = \{ \pm e_i, \,\pm 2e_i,\, 1\leq i \leq q\} \cup\{ \pm e_i\pm e_j, \,1\leq i < j \leq q\}\]
where $(e_1, \ldots, e_q)$ denotes the standard basis of $\mathbb R^q.$
The corresponding Weyl group $W$ is the hyperoctahedral group, which is generated by permutations and sign changes of the $e_i$. 
Put  $R:= \{2\alpha: \alpha\in BC_q\}$ and $\,R_+ := \{ 2e_i, \, 4e_i, \,1\leq i \leq q\} \cup\{ 2(e_i\pm e_j), \,1\leq i < j \leq q\}$  
and denote the associated hypergeometric function  by $F_{BC_q}$. Let $m=m_{p,d}$ be one of the multiplicity functions on $BC_q$ in the geometric cases according to table \eqref{tabelle}  and define $k= k_{p,d}$ on $R$ by
\[ k_{p,d}(2\alpha)  = \frac{1}{2} m_{p,d}(\alpha), \quad \alpha  \in BC_q.\]

Writing $k$ in the form $k=(k_1, k_2, k_3)$ where $k_1$ and $k_2$ are the values on the roots $\pm 2e_i$ and $\pm 4e_i$, respectively and  $k_3$ is the value on the roots $2(\pm e_i\pm e_j)$, we have 
\[ k_{p,d}=\, \bigl(d(p-q)/2, (d-1)/2, d/2\bigr).     \]
The spherical functions of  $G/K$ are then indexed by spectral parameters $\lambda\in \mathbb C^q$ and given by
\[ \varphi_\lambda (a_t) = \widetilde\varphi_\lambda (t) = F_{BC_q}(i\lambda, k_{p,d};t), \quad t\in C.\]
This follows from the fact that for $k=k_{p,d}$,  the commutative algebra
$\{D(p,k); \, p \in S(\mathbb C^q)^W\}$ just represents the radial parts of the algebra
of all invariant differential operators on $G/K$, see Remark 2.3. of \cite{H}.

\begin{example}\textbf{The rank one case.} Here $R_+=\{2, 4\}\subset \mathbb R$. We have multiplicities $k_1, k_2$ and $\rho = \rho(k) = k_1+ 2k_2$. According to the example in \cite{O1}, p.89f, the associated hypergeometric function
is given by
\[ F_{BC_1}(\lambda,k;t) = \, _2F_1\Bigl(\frac{\lambda + \rho}{2}, \frac{-\lambda + \rho}{2}, k_1+k_2+\frac{1}{2}; -\sinh^2 t\Bigr).\]
With $\,\alpha := k_1+k_2-\frac{1}{2}, \, \beta := k_2-\frac{1}{2}$ and the Jacobi functions 
$\phi_\lambda^{(\alpha,\beta)}$ as in \cite{K}, 
this can be written as
\[ F_{BC_1}(i\lambda,k;t) = \phi_{\lambda}^{(\alpha,\beta)}(t).\]
The geometric cases correspond to $\alpha= \frac{dp}{2} -1, \, \beta = \frac{d}{2}-1.$ 
In  Propos.   \ref{prodspher}, the $U_0(1)$-integral cancels (use the  coordinate transform $\widetilde w:= v^{-1}w $),  and the product formula reduces to
\begin{align*}  \widetilde\phi(t)\widetilde\phi(s) 
= &\frac{1}{\kappa_{pd/2}} \int_{B_1} \widetilde\phi\bigl({\rm arcosh} |\!\cosh t \cosh s   \,+\, w\sinh t \sinh s|)\cdot (1-|w|^2)^{\frac{pd}{2}-\gamma} dw\\
= & \int_{\Sigma_{p,1}}\! \widetilde\phi\bigl({\rm arcosh} |\!\cosh t \cosh s   + x_1\sinh t \sinh s|)d\sigma(x)\end{align*}
where $\widetilde\varphi = \varphi_\lambda^{(\alpha, \beta)}$ with  $\alpha =  \frac{pd}{2} -1, \, \beta = \frac{d}{2}-1.$ The second identity is obtained by formula \eqref{splitting} 
for the sphere $\Sigma_{p,1}= \{x\in \mathbb F^p: |x| = 1\}.$ 
In view of relation (5.24) in \cite{K}, this  formula just coincides with the product formula in rank 1 given in Section 7 of \cite{K}, 
\begin{align}\label{FJK} \phi_\lambda^{(\alpha, \beta)}(t)\phi_\lambda^{(\alpha, \beta)}(s) \, = \,c_{\alpha,\beta}
\int_0^1 \int_0^\pi   &\phi_\lambda^{(\alpha, \beta)}\bigl({\rm arcosh} |\!\cosh t \cosh s   \,+\, re^{i\psi}\sinh t \sinh s|)\cdot\notag \\
\cdot& (1-r^2)^{\alpha-\beta-1}r^{2\beta + 1} (\sin\psi)^{2\beta} 
r dr d\psi\end{align}
which degenerates for $\beta=-1/2$ (i.e. $\mathbb F=\mathbb R$) to an integral over $[-1,1]$ with respect to  $\,(1-r^2)^{\alpha-1/2} dr.$ 

In fact,  formula  \eqref{FJK} was established in \cite{FK}   for arbitrary $\alpha\geq\beta\geq -\frac{1}{2}$ with $ (\alpha, \beta) \not=\bigl(-\frac{1}{2}, -\frac{1}{2}\bigr),$ i.e. arbitrary non-negative root multiplicities different from zero.

\end{example}

\section{Continuation of the product formula}\label{prod_formula}

In the following, $q$ and $d=\text{dim}_{\mathbb R}(\mathbb F)$ are fixed. 
For $\mu \in \mathbb C$ with $\text{Re}\, \mu>\gamma -1$ and spectral parameter $\lambda\in \mathbb C^q$
define 
\[ \varphi_\lambda^\mu(t) = F_{BC_q}(i\lambda, k_\mu;t);  \quad t\in \mathbb R^q \]
with multiplicity
\[ k_\mu = \bigl(\mu-dq/2, (d-1)/2, d/2\bigr).\]
If $\mu= pd/2$, then $k_\mu = k_{p,d}$ as in the previous section.

\begin{theorem}\label{prodformel} For $\mu\in \mathbb C$ with $\text{Re}\,\mu>\gamma -1$, the hypergeometric functions $\varphi_\lambda^\mu$ satisfy the product formula
\[\varphi_\lambda^\mu(t)\varphi_\lambda^\mu(s)\,=\, (\delta_t*_\mu\delta_s)(\varphi_\lambda^\mu) \]
with the probability measures
\[( \delta_t*_\mu\delta_s)(f) = 
\frac{1}{\kappa_{\mu}} \int_{B_q}\int_{U_0(q,\mathbb F)} f\bigl(d(t,s;v,w)\bigr) \Delta(I-w^*w)^{\mu-\gamma} dv dw \]
where $\kappa_\mu$ is given by \eqref{kappa} and the argument is
\[ d(t,s;v,w) = {\rm arcosh}({\rm spec}_s(\sinh\underline t\, w \sinh\underline s + \cosh\underline t \,v \cosh\underline s)).\]
\end{theorem}

This is a partial generalization of formula \eqref{FJK} by Flensted-Jensen and Koornwinder   for $BC_1$ to higher rank.

\begin{proof} The basic idea is analytic continuation with respect to the parameter $\mu$ in the right half plane by use of 

\smallskip\noindent
\textbf{Carlson's Theorem} \,(see e.g. ([T], p.186):
Let $f$ be a function which is holomorphic in a neighbourhood of $\{ z\in \mathbb C: \text{Re}\, z \geq 0\}$ 
satisfying $\, f(z) = O(e^{c|z|})$ for some constant $c<\pi$. Suppose that $f(n) = 0$ for all $n\in \mathbb N_0$. Then $f$ is identically zero.
\smallskip

\smallskip\noindent
A direct application of Carlson's theorem would require moderate exponential growth of the hypergeometric function with respect to the relevant multiplicity 
parameter $k_1$ in a right half plane. So far however, sufficient exponential estimates are available only for real,  non-negative  multiplicities 
(Propos. 6.1. of \cite{O1}, and the results of \cite{S}). We therefore proceed in two steps. First, we restrict
to a discrete set of spectral parameters, for which the 
hypergeometric function is a Jacobi polynomial and the required growth properties are
guaranteed. In a second step,
we fix a non-negative multiplicity and carry out analytic continuation with respect to the spectral parameter,
using known bounds on the hypergeometric function for non-negative multiplicities.

\smallskip
To go into detail, let  $R^\vee = \{\alpha^\vee: \alpha\in R\}$ be the root system dual to $R,$  $Q^\vee = \mathbb Z.R^\vee$ the coroot lattice and 
$P=\{ \lambda\in \mathbb R^q: \langle\lambda,\alpha^\vee\rangle \in \mathbb Z \,\,\forall\,\alpha\in R\}\,$ the weight lattice of  $R$. Further, denote by
$\, P_+ =\{ \lambda\in P: \langle\lambda,\alpha^\vee\rangle \geq 0 \,\,\forall\,\alpha\in R_+\}\,$  the set of dominant weights associated with $R_+$.
Then for $k\in K^{reg}$ and $\lambda\in P_+$,
\[ F_{BC_q}(\lambda +\rho(k), k;t) = \,c(\lambda+\rho(k),k) P_\lambda(k;t)\]
where $c(\lambda,k)$ is the $c$-function \eqref{c_function} which is meromorphic on $\mathbb C^q\times K$, and the 
 $P_\lambda$ are the Heckman-Opdam Jacobi polynomials of type $BC_q$; see \cite{HS}, equation (4.4.10).
In our case,  $\rho(k)$ is given by 
\[ \rho(k) = (k_1 + 2k_2)\sum_{i=1}^q e_i + 2k_3\sum_{i=1}^q(q-i)e_i = (\mu-\frac{dq}{2} +d-1)\sum_{i=1}^q e_i\,+ d \sum_{i=1}^q (q-i)e_i.\]
Using the asymptotics of the gamma function, one checks that for fixed $\lambda\in P_+$, the function
$\, c(\lambda+\rho(k_\mu),k_\mu)\,$ is bounded away from zero 
as $\mu\to\infty$ in the right half plane 
\[H=\{\mu\in \mathbb C: \text{Re}\,\mu>\gamma-1\}.\] 
Indeed, for $\rho = \rho(k)$ with $k= (k_1,k_2, k_3)$ one has
\begin{align*}
 &c(\lambda + \rho,k) \,=\\
&=\prod_{i=1}^q \frac{\Gamma(\lambda_i + \rho_i)\,\Gamma(\rho_i + k_1)}
{\Gamma(\lambda_i + \rho_i +k_1)\,\Gamma(\rho_i)} \,\cdot 
\prod_{i=1}^q \frac{\Gamma\bigl(\frac{\lambda_i+\rho_i}{2} + \frac{1}{2}k_1\bigr)\,
\Gamma(\frac{\rho_i}{2} + \frac{1}{2} k_1 + k_2)}{\Gamma\bigl(\frac{\lambda_i+\rho_i}{2} + \frac{1}{2}k_1+k_2\bigr)\,\Gamma(\frac{\rho_i}{2} + \frac{1}{2} k_1)}\\
&\cdot\prod_{i<j} \frac{\Gamma\bigl(\frac{\lambda_i + \rho_i - \lambda_j -\rho_j}{2}\bigr)\,
\Gamma\bigl(\frac{\rho_i - \rho_j}{2} + k_3\bigr)}
{\Gamma\bigl(\frac{\lambda_i + \rho_i - \lambda_j - \rho_j}{2} + k_3\bigr)\, \Gamma\bigl(\frac{\rho_i-\rho_j}{2}\bigr)}\,\cdot
\,\prod_{i<j}\frac{\Gamma\bigl(\frac{\lambda_i + \rho_i + \lambda_j +\rho_j}{2}\bigr)\,
\Gamma\bigl(\frac{\rho_i + \rho_j}{2} + k_3\bigr)}
{\Gamma\bigl(\frac{\lambda_i + \rho_i +\lambda_j + \rho_j}{2} + k_3\bigr) \, \Gamma\bigl(\frac{\rho_i+\rho_j}{2}\bigr)}.
\end{align*}
As $k_1\to \infty$ in the half-plane $\, \text{Re}\, k_1 > 0$,  the first product is asymptotically equal to $\prod_{i=1}^q 
\bigl(\frac{1}{2}\bigr)^{\lambda_i}$, the second one is asymptotically equal to $1$, the third product is independent of $k_1$, and the last product is asymptotically equal to $1$. Thus for fixed $\lambda, \,\, c(\lambda + \rho,k)$ is bounded away from zero.

\smallskip

According to Proposition \ref{prodspher}, the $P_\lambda(k_\mu;.)$ with $\mu=pd/2\, (p\geq 2q)$  satisfy the product formula \begin{align}\label{prodpoly}
&P_\lambda(k_\mu;t)P_\lambda(k_\mu;s)\notag \\
\, \,\,& =\frac{1}{\kappa_{\mu}} \int_{B_q}\int_{U_0(q,\mathbb F)} \frac{1}{c(\lambda+\rho(k_\mu),k_\mu)}P_\lambda\bigl(k_\mu;d(t,s;v,w)\bigr) \Delta(I-w^*w)^{\mu-\gamma} dv dw\end{align}
for all $t,s\in \mathbb R^q.$ 
The Jacobi polynomials $P_\lambda(k;.)$ have rational coefficients in $k$ with respect to the monomial  basis $e^\nu,\,\nu\in P$. This is shown in Par.11 of \cite{M}, but it also follows from the explicit determinantal construction in \cite{DLM}, Theorem 5.4.
Moreover, as derived in the proof of Theorem 3.6 of \cite{R}, the normalized integral
\[ \frac{1}{|\kappa_{\mu}|} \int_{B_q} \vert\Delta(I-w^*w)^{\mu-\gamma}\vert dw\]
converges exactly if $\,\text{Re}\, \mu >\gamma -1$ and 
is  of polynomial growth as $\mu\to\infty$ in $H$.
Thus for fixed $t,s,$ both sides of \eqref{prodpoly} are holomorphic in $\mu\in H$ and of 
polynomial growth as $\mu\to\infty$ in $H$. Moreover, they coincide for all half integer values $\mu = pd/2, \, p\geq 2q$. Application of Carlson's theorem yields that formula \eqref{prodpoly} holds for all $\mu\in H$. 
This proves the stated result for spectral parameters $\lambda + \rho(k)$ with $\lambda\in P_+$ and $k=k_\mu, \, \mu\in H$. 

\smallskip 
Denote again by $C\subset \mathbb R^q$  the closed Weyl chamber associated with $R_+$.
In order to extend the product formula with respect to the spectral parameter, we fix $s,t\in C$ as well as $k=k_\mu$ and restrict to 
real $\mu>\gamma-1$ first. Then $k$ is nonnegative, and we have the following exponential estimate for $F_{BC_q}$ from \cite{O1}, Proposition 6.1:
\[ |F_{BC_q}(\lambda,k;t)|\,\leq \,|W|^{1/2} e^{{max}_{w\in W}\text{Re}\langle w\lambda,t\rangle} .\]
Let $\,H^\prime := \{\lambda\in \mathbb C^q: \text{Re}\lambda\in C^0\}\,$. Then for $\lambda\in H^\prime$ and all $w\in W$, 
\[\text{Re}\langle w\lambda,t\rangle \,\leq \text{Re}\langle\lambda,t\rangle.\]
Choose a constant vector $a\in C^0$ so large that $d(t,s;v,w) -a\,$ is contained in the negative chamber $-C\,$ 
for all $v\in U(q)$ and all $w\in B_q$. Then consider
\[\widetilde F(\lambda,k;t):= \, e^{-\langle\lambda, a+t\rangle} F_{BC_q}(\lambda,k;t).\]
The function $\widetilde F$ is bounded as a function of  $\lambda\in H^\prime$. 
If the spectral parameter is of the form $\lambda = \widetilde \lambda + \rho(k_\mu)$ with $\widetilde\lambda\in P_+$, then by our previous results we have the product formula
\begin{align}\label{tildeformel}
& \widetilde F(\lambda,k_\mu;t)\widetilde F(\lambda,k_\mu;s) \,=\notag \\
&=\frac{1}{\kappa_{\mu}} \int_{B_q}\int_{U_0(q,\mathbb F)}\!
e^{\langle\lambda,\, d(t,s;v,w)-a-s-t\rangle} \widetilde F(\lambda,k_\mu;d(t,s;v,w))\Delta(I-w^*w)^{\mu-\gamma} dv dw.\end{align}
Both sides are holomorphic and bounded in $\lambda\in H^\prime$. 
We are now going to carry out analytic extension with respect to $\lambda$. For this, 
we choose a set of fundamental weights  $\{\lambda_1, \ldots, \lambda_q\}\subset P_+$ and write $\lambda\in H^\prime$ as 
$\lambda = \sum_{i=1}^q z_i\lambda_i$ with coefficients $z_i\in \{z\in \mathbb C: \text{Re} z >0\}$. 
Successive holomorphic extension with respect to $z_1, \ldots , z_q$ by use of Carlson's Theorem then
yields the validity of \eqref{tildeformel} for all $\lambda \in H^\prime$, and thus, by $W$-invariance and continuity, for all $\lambda\in \mathbb C^q$. This proves  the stated product formula for real $\mu > \gamma -1$. Analytic continuation finally gives it for all $\mu\in H$, which finishes the proof of Theorem \ref{prodformel}.

\end{proof}

\section{Hypergroup algebras associated with $F_{BC}$}\label{hypergroups}

The positive product formula of Theorem \ref{prodformel} leads to three continuous series ($d=1,2,4$) of positivity-preserving convolution algebras on the Weyl chamber $C$ which are parametrized by $\mu$. We shall describe them as commutative hypergroups, having Heckman-Opdam hypergeometric functions as characters. 
In the group cases, which correspond to the discrete values $\mu=pd/2$, these hypergroup algebras
are just given by the double coset convolutions associated with the Gelfand pairs $(G,K)$ as in Section \ref{spherical}. 
In the rank one case, they coincide with the well-known one-variable Jacobi hypergroups.

Let us first briefly recall some key notions and facts from hypergroup theory. For a detailed treatment, the reader is referred to \cite{J}. 
Hypergroups generalize the convolution algebras of locally compact groups, 
with the convolution product of two point measures $\delta_x$ and $\delta_y$ being in general not a point measure again
but  a probability measure depending on $x$ and $y$.

\begin{definition}
A hypergroup is a locally compact Hausdorff space $X$ with a weakly continuous, associative convolution $*$ on the space
 $M_b(X)$ of regular bounded Borel measures on $X$, satisfying the following properties:
\begin{enumerate}\itemsep=-1pt
\item 
The convolution product $\delta_x*\delta_y$ of two point measures is a compactly supported probability measure on $X$, and  $\text{supp}(\delta_x*\delta_y)$ depends continuously on $x$ and $y$ with respect to the so-called Michael topology on the space of compact subsets of $X$ (see \cite{J}).
\item There is a neutral element $\delta_e$ satisfying $\delta_e*\delta_x = \delta_x = \delta_x*\delta_e$ for all $x\in X$. 
\item There is 
  a continuous involution $x\mapsto\bar x$ on $X$ such that for all $x, y\in X$,
$e\in \text{supp}(\delta_x*\delta_y)$ is equivalent to $ x=\bar y$,
and
$\delta_{\bar x} * \delta_{\bar y} =(\delta_y*\delta_x)^-$.
Here for $\mu\in M_b(X)$, the measure  $\mu^-$ is given by
$\mu^-(A)=\mu(A^-)$ for Borel sets $A\subset X$.
\end{enumerate}
\end{definition}

Due to weak continuity, the convolution of measures on a hypergroup is uniquely determined by the convolution of point measures. 

\smallskip
If the convolution is commutative, then
$(M_b(X),*)$ becomes a commutative Banach-$*$-algebra with identity $\delta_e$.  
Moreover, there exists an (up to a multiplicative factor) unique Haar measure $\omega$, that is a positive Radon measure on $X$ satisfying 
\[ \int_X f(x*y)d\omega(y)  = \int_X f(y)d\omega(y)  \quad \text{for all } \, x\in X, \, f\in C_c(X),\]
where $f(x*y)= (\delta_x*\delta_y)(f).\,$  
The multiplicative functions of a commutative hypergroup $X$ are given by
\[ \chi(X) = \{\varphi\in C(X): \,\varphi\not= 0, \, \varphi( x* y ) = \varphi(x)\varphi(y) \,\, \forall\,  x,y\in X\}.\]  
The decisive object for harmonic analysis is the dual space of $X$, 
defined by
\[ \widehat X :=\{\varphi\in \chi(X): \,\varphi \text{ is bounded and }\, \varphi(\overline x)= \overline{\varphi(x)}\,\, \forall\,  x\in X\}.\] 
The elements of $\widehat X$ are called characters. As in the case of LCA groups, the dual of a commutative hypergroup is a locally compact Hausdorff space with the topology of locally uniform 
convergence. It is naturally identified with the symmetric part of the spectrum of the convolution algebra $L^1(X,\omega)$. In contrast to the group case, $\widehat X$ is often a proper subset of $\chi(X).$ 
The  Fourier transform on $ L^1(X, \omega)$ is defined by $\,\widehat f(\varphi) := \int_X f\overline{\varphi}d\omega.\,$
It is injective, and there exists a 
unique positive Radon measure $\pi$ on $\widehat X$, called the Plancherel measure of $(X,*)$, such that $f\mapsto \widehat f$ extends to an isometric 
isomorphism from $L^2(X,\omega)$ onto $L^2(\widehat X,\pi)$. 
 As for groups, there are convolutions between functions from various classes of $L^p$-spaces (or measures) on a hypergroup with Haar measure $\omega$. For example, if $1\leq p\leq \infty$ and $f\in L^1(X,\omega), g\in L^p(X,\omega)$, then the convolution product
\[ f*g(x) = \int_X f(x*\overline y) g(y) d\omega(y)\]
belongs to $L^p(X,\omega)$ and satisfies $\, \|f*g\|_{p,\omega}\leq \|f\|_{1,\omega}\|g\|_{p,\omega}.$ 

\medskip
Let us come back to the situation of Section 3. With the notions from there, we can now state our main theorem:

\begin{theorem} \begin{enumerate}\itemsep=-1pt
\item[\rm{(1)}] Let $\mu > \gamma -1$. Then the probability measures given by
\begin{equation}\label{convoformel} (\delta_s *_\mu \delta_t)(f) = \frac{1}{\kappa_{\mu}} \int_{B_q}\int_{U_0(q,\mathbb F)} f\bigl(d(s,t;v,w)
\bigr) \Delta(I-w^*w)^{\mu-\gamma} dv dw \end{equation}
for $\,s,t \in C$ define a commutative hypergroup structure $C_\mu=(C,*_\mu)$  on the chamber $C \cong \overline{\frak a_+}.$ The neutral element is $0$ and  the involution is the identity mapping. The support of $\delta_s*_\mu\delta_t$ satisfies
\[ \text{supp}(\delta_s*_\mu\delta_t) \,\subseteq \{r\in C: \|r\|_\infty \leq \|s\|_\infty + \|t\|_\infty\}\]
where $\|\,.\|_\infty$ is the maximum norm in $\mathbb R^q$. 
\item[\rm{(2)}]  A Haar measure of the hypergroup $C_\mu$ is given by the weight function \eqref{weight} of the corresponding hypergeometric transform,
\begin{align*} d\omega_\mu(t) \,=\,  const\,\cdot&\prod_{i=1}^q |\sinh t_i|^{2\mu-d(q-1)-1} |\cosh t_i|^{d-1}\,\cdot \\
\cdot\! &\prod_{1\leq i <j\leq q} \! |\cosh(2t_i) - \cosh(2t_j) |^d\,dt. \end{align*}
\end{enumerate}
\end{theorem}

\begin{proof} (1)
It is clear that $\delta_s*_\mu\delta_t$ is a probability measure on $C$ with 
\begin{align*} \text{supp}(\delta_s*_\mu\delta_t) = \,\{ d(s,t;v,w)=\text{arcosh}
({\rm spec}_s(\sinh\underline s\, w &\sinh\underline t + \cosh\underline s \,v \cosh\underline t)),\\
& v\in U_0(q,\mathbb F), w\in \overline{B_q}\}.\end{align*}
For the support statement, we 
denote by $\|A\|$ the spectral norm of $A\in \mathbb F^{q\times q}$, that is $\, \|A\| = \|{\rm spec}_s(A)\|_\infty$ (the biggest singular value of $A$).
By the submultiplicativity of $\|\,.\,\|$ we obtain for $v$ and $w$  within the relevant range the estimate
\begin{align*} 
\|\sinh\underline s\, w \sinh\underline t + \cosh\underline s \,v \cosh\underline t\| \,\leq &\,
\|\sinh\underline s\|\|\sinh\underline t\| + \|\cosh\underline s\| \|\cosh\underline t\| \\
 = &\, \sinh\!\|s\|_\infty \cdot \sinh\!\|t\|_\infty\, + \,\cosh\!\|s\|_\infty\cdot \cosh\!\|t\|_\infty\\
=& \cosh(\|s\|_\infty + \|t\|_\infty).\end{align*}
This implies the stated support inclusion. For the weak continuity of the convolution $*_\mu$ on $M_b(C)$, it 
suffices to verify that for each $f\in C_b(C),$ the mapping $\, (s,t)\mapsto f(s*_\mu t)$ is continuous. 
But this is immediate because $d(s,t;v,w)$ depends continuously on its arguments. To see that $*_\mu$ is commutative, we note that  ${\rm spec}_s(A) = 
{\rm spec}_s(A^*)$ for $A\in \mathbb F^{q\times q}$, and hence $d(t,s;v,w) = d(s,t;v^*,w^*)$. As the integral in \eqref{convoformel} is invariant under the substitution $v\mapsto v^* = v^{-1}, w\mapsto w^*$, it follows that $\,\delta_t*_\mu\delta_s = \delta_s*_\mu\delta_t$. 
For the associativity of $*_\mu$ it suffices to verify that 
\[\delta_r*_\mu(\delta_s*_\mu\delta_t)(f) = (\delta_r*_\mu\delta_s)*_\mu\delta_t(f)\]
for all $f\in C_c^\infty(\mathbb R^q)^W$ and all $r,s,t\in C$. 
In view of the Paley-Wiener theorem for the hypergeometric transform, both sides are equal to
\[\int_{i\mathbb R^q} \mathcal Ff(\lambda) F(\lambda,k;r) F(\lambda,k;s) F(\lambda,k;t)d\nu(\lambda).\]
This proves the assertion.

From the explicit form of the convolution it is obvious that $0$ is neutral. 
In the discrete cases $\mu = pd/2$ coming from Gelfand pairs, $*_\mu$ is the convolution of a double coset hypergroup. Moreover, $\text{supp}(\delta_s*_\mu\delta_t)$  is independent of $\mu$. In order to see that the identity mapping is a hypergroup  involution for all $\mu$,   it therefore suffices (by uniqueness of an  involution) to show that the zero matrix $0$ is contained in $\text{supp}(\delta_t*_\mu\delta_t).$ But
\[ d(t,t; I_q, -I_q) = {\rm arcosh}({\rm spec}_s(-(\sinh\underline t)^2\,+ (\cosh\underline t)^2)) \, = \, {\rm arcosh} (I_q) = 0,\]
which proves the claim.

\smallskip\noindent
(2) Let $f, g\in C_c^\infty(\mathbb R^q)^W$. Notice first that
\[ f(t) = \int_{i\mathbb R^q} \mathcal F f(\lambda)F(\lambda,k; t) d\nu(\lambda) \]
by the inversion theorem for the hypergeometric transform (Theorem \ref{Plancherel}).
As $\, F(\lambda,k; s*_\mu t) = F(\lambda,k;s) F(\lambda,k;t)\,\,$ for all $s,t\in C\,$ we obtain, with the notation of 
Corollary \ref{translate}, 
\[ f(s*_\mu t) = \tau_sf(t). \]
The Plancherel formula (Theorem \ref{Plancherel}) further gives
\begin{align*} \int_C (\tau_s f)\overline g\, d\omega_\mu\,= &\, c\int_{iC} \mathcal F(\tau_s f) \overline{\mathcal F g}\, d\nu\,=\, c\int_{iC} \mathcal F f(\lambda)F(\lambda,k;s) \overline{\mathcal F g(\lambda)}\,d\nu(\lambda) \\
=\,& c \int_{iC} \mathcal F f(\lambda)\,\overline{\mathcal F(\tau_s g)(\lambda)}\,d\nu(\lambda)\,=\,
\int_C f(\overline{\tau_sg}) d\omega_\mu
\end{align*}
with a constant $c>0$. It was used here that $\overline{F(\lambda,k;s)} = F(-\lambda,k;s) = F(\lambda,k;s)$ for $\lambda\in iC.$ 
Choose now a sequence $g_n \in C_c^\infty(\mathbb R^q)^W, \, n\in \mathbb N$ such that $g_n\uparrow 1$ pointwise. Then also $\,\tau_s(g_n) \uparrow 1$, and the monotonic convergence theorem shows that
\[ \int_{C} (\tau_s f) d\omega_\mu\,=\, \int_{C} f d\omega_\mu.\]
This proves that $\omega_\mu$ is a Haar measure of $C_\mu$. 

\end{proof}

\begin{lemma} Suppose that $\varphi: C_\mu\to \mathbb C$ is continuous and multiplicative, i.e.
\[\varphi(s)\varphi(t) = \varphi(s*_\mu t) \quad \text{for all } s,t \in C. \]
Then $\,\varphi = \varphi_\lambda^\mu$ with some  $\lambda\in \mathbb C^q$. 
\end{lemma}

\begin{proof} The proof follows standard arguments.
For abbreviation, we write $k=k_\mu,\, \omega=\omega_\mu$ and $\,*=*_\mu.$  In a first step, consider $g\in C_c^\infty(\mathbb R^q)^W.$ Let $p\in S(\mathbb C^q)^W$ be a $W$-invariant polynomial and $T(p) = T(p,k_\mu)$ the associated
Cherednik operator.
As 
\[ g(s*t) = \int_{i\mathbb R^q} \mathcal Fg(\lambda) F(\lambda,k; s) F(\lambda, k; t) d\nu(\lambda)\]
for all $s,t\in C$, we obtain
 \begin{equation}\label{convodiff}
T(p)_s\, g(s*t) = \int_{i\mathbb R^q} \mathcal Fg(\lambda) p(\lambda) F(\lambda,k; s) F(\lambda, k; t) d\nu(\lambda)\,=\, 
T(p)_t\, g(s*t)\end{equation}
and 
\[T(p)_s\,g(s*t)\vert_{s=0} \,=\, T(p)g(t).\]
Suppose now $\varphi $ is continuous, non-zero and multiplicative on $C$. Notice first that $\varphi(0)=1$, because $0$ is neutral. We extend $\varphi$ to a $W$-invariant function on $\mathbb R^q$ and choose $g\in C_c^\infty(\mathbb R^q)$ with $\,\int_C \varphi g\,d\omega\, = 1.$
Recall that the involution of the hypergroup $C_\mu$ is the identity. Thus
\[ \varphi *g (s) =\, \int_C \varphi(s*t)g(t)d\omega(t)\,=\, \varphi(s) \]
and therefore 
\[\varphi(s) = \varphi*g(s) = \int_C \varphi(t) \tau_tg(s)d\omega(t),\] which belongs to $C^\infty(\mathbb R^q)$ because $\tau_tg\in C_c^\infty(\mathbb R^q)$ for all $t$ according to Lemma \ref{translate}. Further, 
\[ \varphi(s*t) = \int_C \varphi(r)\tau_rg(s*t)\,d\omega(r)\]
and therefore
\begin{align*} T(p)_s\,\varphi(s*t) &= \int_C \varphi(r) T(p)_s (\tau_r g(s*t)) d\omega(r) \,=\, \int_C \varphi(r) \,T(p)_t (\tau_r g(s*t)) \,d\omega(r)\\
 &= T(p)_t\, \varphi(s*t).
\end{align*}
In particular,
\[ T(p)\varphi(t)\,=\, T(p)_t\, \varphi(s*t)\vert_{s=0}\,=\,T(p)_s\, \varphi(s*t)\vert_{s=0}\, = \sigma_\varphi(p)\cdot \varphi(t)\]
with  $\, \sigma_\varphi(p) = (T(p)\varphi)(0).$ 
The mapping $p\mapsto \sigma_\varphi(p)$ is obviously multiplicative and linear on $S(\mathbb C^q)^W$. According
to a well-known result form invariant theory (see e.g. \cite{Hel}, Ch. III.4, Lemma 3.11), it coincides with a point evaluation, that is 
\[ \exists \,\lambda\in \mathbb C^q: \, \sigma_\varphi(p) = p(\lambda) \,\,\,\forall \,p \in S(\mathbb C^q)^W.\]
It is thus shown that $\varphi$ satisfies the hypergeometric system \eqref{hypersystem} withspectral parameter $\lambda$, corresponding to $R=BC_q$ and $k=k_\mu$. By uniqueness of the solution, it follows that $\, \varphi= F_{BC_q}(\lambda,k_\mu; \,.) =  \varphi_{-i\lambda}^\mu.$

\end{proof}

\begin{theorem} \label{main_Theorem} The set of multiplicative functions and the dual space of the hypergroup $C_\mu$ are given by
\begin{align*} \chi(C_\mu) = &\,\{\varphi_\lambda = \varphi_\lambda^\mu: \,\lambda\in C+iC\};\\ \notag
\widehat{ C_\mu} =& \,\{ \varphi_\lambda\in \chi(C_\mu): \, \overline\lambda\in W.\lambda \,\,\text{ and }\,  \, {\rm Im} \,\lambda \in \rm{co}(W.\rho)\}
\end{align*}
where $\rho = \rho(k_\mu)$ and $ {\rm co}(W.\rho)$ denotes the convex hull of the Weyl group orbit $W.\rho$.
\end{theorem}

The second part of this theorem is in accordance with the characterization of the bounded spherical functions of a Riemannian symmetric space of non-compact type, see \cite{Hel}, Chap. IV,  Theorem 8.1. In our more general context, we shall not work with an integral representation
but proceed by using estimates on the hypergeometric function given in \cite{S} as well as the generalized Harish-Chandra expansion of \cite{O1}. We mention at this point that for the Grassmann manifolds over $\mathbb F= \mathbb R$, there is an explicit integral formula for the spherical functions given in \cite{Sa} which could probably also be used after analytic extension.

\begin{proof} (Proof of Theorem \ref{main_Theorem}.)
The identification of $\chi(C_\mu)$ is furnished by the previous lemma.
For the identification of the dual space, note first that 
\[\overline{\varphi_\lambda} = \varphi_{\overline\lambda}\]
as a consequence of \eqref{realF}. Thus $\varphi_\lambda$ is real if and only if $\, \overline \lambda\in W.\lambda$.
It remains to identify those functions from $\chi(C_\mu)$ which are bounded. For this, we
observe first that the set $A=\{\lambda \in \mathbb C^q: \varphi_\lambda\in \widehat{C_\mu}\}$ is closed in $\mathbb C^q$. Indeed, suppose that $(\lambda_i)_{i\in \mathbb N} \,$ is a sequence in $A$ which converges to $\lambda_0\in \mathbb  C^q$. 
Being members of a hypergroup dual, the $\varphi_{\lambda_i}$ are uniformly bounded by $1$ (see \cite{J}).
As $\, (\lambda, t)\mapsto \varphi_\lambda(t)$ is continuous, it follows by a standard compactness argument that the sequence $\varphi_{\lambda_i}$ converges to $\varphi_{\lambda_0}$ locally uniformly on $C$ (see e.g. \cite{D}, Chap.XII, Sec.8). 
This implies  that  $\varphi_{\lambda_0}$  belongs to $\widehat{C_\mu}$ as well.

 \smallskip
\noindent
Next recall that \[\varphi_\lambda(t) = F_{BC_q}(i\lambda,k_\mu;t)=: F_{i\lambda}(t)\]
 and notice that $F_{-\lambda} = F_\lambda$. We thus have to prove that $F_\lambda$ is bounded 
if and only if $\,{\text Re}\,\lambda\in \rm{co}(W.\rho).$
We may assume that  $\lambda= \xi + i\eta$ with $\,\xi,\eta\in C$. By Cor. 3.1  of \cite{S}, 
\begin{equation}\label{eq1} |F_\lambda(t)| \, \leq \, F_{\xi}(t) \quad\forall \,  t\in C.\end{equation}
Further, according to Remark 3.1 of [loc.cit],  $ F_{\xi} $ behaves asymptotically (for large arguments in $C$) as 
\begin{equation}\label{eq2} F_ {\xi}(t) \asymp \, e^{\langle \xi-\rho,t\rangle}\,\cdot\!\!\prod_{\alpha\in R_0^+\vert \langle\alpha,\xi\rangle =0}\!\!
 \bigl(1 + \langle \alpha,t\rangle \bigr).\end{equation}
Here $R_0^+$ are the indivisible positive roots, in our case $\,R_0^+ = \{2e_i, \, 2(e_i\pm e_j), 1\leq i <j\leq q\}.\,$
Consider now  $\lambda = \xi + i\eta$ with $\,\xi={\rm Re} \,\lambda \in \rm{co}(W.\rho).$ We claim that $F_\lambda$ is bounded. By closedness of $A$, it 
suffices to assume that $\,\xi$ is actually contained in the open interior of $\rm{co}(W.\rho).$ Then there exists a constant $0<s<1, s= 1-\epsilon$, such that $\xi\in \rm{co}(W.s\rho).$ We use the characterization 
\begin{equation}\label{eq3}  \text{co}(W.x) = \bigcap_{w\in W} w(x-C^*)\end{equation}
for $x\in C$, where $\,C^* = \{ x\in \mathbb R^q: \langle t,x\rangle \geq 0 \,\,\forall t\in C\}$ is the closed dual cone of $C$; see e.g. \cite{Hel}, Lemma IV.8.3. This shows that $\, s\rho-\xi\in C^*$ and therefore 
\[ \langle \xi -\rho, t\rangle = \langle \xi -s\rho ,t \rangle - \epsilon\langle \rho,t\rangle \,\leq \, - \epsilon\langle \rho,t\rangle \quad \forall \, t\in C.\]
Note that  $\langle \rho,t\rangle >0$  
for all $t\in C\setminus\{0\}$, because  our multiplicity is non-negative and different from zero. Hence $ \langle \rho,t\rangle \geq c|t|$ for some constant $c>0$. Together with estimates \eqref{eq1} and \eqref{eq2}, this proves  boundedness of $F_\lambda$ as claimed.

\smallskip\noindent
For the converse inclusion, we have to show that $F_\lambda$ is unbounded if $\xi = \text{Re}\lambda \notin \rm{co}(W.\rho).$ For 
real $\lambda=\xi\in C$ we use again  \eqref{eq2}. According to \eqref{eq3}, 
there exists some $t\in C$ such that $\,\langle \xi-\rho, t\rangle >0$ (recall that $\langle \xi,t\rangle \geq \langle w\xi,t\rangle$ for all 
$w\in W$). This implies that $F_\xi$ is unbounded in $C$. 

\smallskip\noindent In case $\eta= \text{Im}\lambda\not=0$ we employ
the Harish-Chandra expansion of
$F_\lambda$ (see \cite{O1}) in the interior $C^\circ$ of $C$. It is of the form
\[ F_\lambda(t) = \sum_{w\in W} c(w\lambda) e^{\langle w\lambda -\rho,t\rangle} \Bigl(\sum_{q\in Q^+} \Gamma_q(w\lambda) e^{-\langle q,t\rangle}\Bigr)\]
with (unique) coefficients $\Gamma_q(w\lambda)\in \mathbb C$, where $\Gamma_0(w\lambda) =1.$  Here $Q^+$ is the positive lattice generated by  $R_+$ and $c(\lambda) = c(\lambda,k_\mu)$ denotes the $c$-function. 

As $\xi\in C\setminus  \rm{co}(W.\rho),$ there is some $t\in C$ and hence also some $t_0\in C^\circ$ such that $\langle \xi-\rho,t_0\rangle >0.$ 
Fix  $t_0$ and consider $F_\lambda(st_0)$ for $s\in \mathbb R, s\to +\infty$.  As the imaginary part of $\lambda$ is nonzero, Lemma 4.2.2. in Part I of \cite{HS} implies that there exist constants
$M_{w\lambda }>0$ (depending on $t_0$) such that 
\[ |\Gamma_q(w\lambda)|\leq M_{w\lambda}e^{\langle q, t_0\rangle} \quad \text{ for all }\, q\in Q^+.\]
For $s\in \mathbb R, s>0$ we may therefore estimate
\[ \big\vert \sum_{q\in Q^+\setminus\{0\}} \Gamma_q(w\lambda) e^{-\langle q, st_0\rangle}\big\vert \,\leq M_{w\lambda} \sum_{q\in Q^+\setminus\{0\}} e^{(1-s)\langle q,t_0\rangle},\]
which tends to zero as $s\to +\infty$. 
 Thus
\[ F_\lambda(st_0)  \asymp \sum_{w\in W} c(w\lambda) e^{\langle w\lambda -\rho,st_0\rangle} \,=\, 
\sum_{w\in W} c(w\lambda) e^{s\langle w\xi -\rho,t_0\rangle} e^{is\langle w\eta,t_0\rangle}
\quad \text{as }\, s\to +\infty. \]
Notice that  $c(\lambda) \not=0.$ Moreover,  $\,\langle \xi-\rho,t_0\rangle \geq
\langle w\xi-\rho,t_0\rangle $ for all $w\in W$ where equality can only occur if $w\xi = \xi.$ Therefore,  the leading term of the last sum   is
\[ e^{s\langle \xi-\rho,t_0\rangle}\cdot \sum_{w\in W_\xi} c(w\lambda) e^{is\langle w\eta,t_0\rangle }\]
with  $\, W_\xi = \{ w\in W: w\xi = \xi\}.$ Application of Lemma \ref{MV} below now implies that     $ s\mapsto F_\lambda(st_0)$ is   unbounded as $s\to +\infty$. This finishes the proof.
\end{proof}

\begin{lemma}\label{MV}
 Let 
$\,  f(s) = e^{as} \cdot \sum_{k=1}^N c_k\, e^{i\lambda_k s} \quad \text{on }\,  \mathbb R\,$
with constants $a>0,$ $c_k\in \mathbb C$ which are not all zero, and distinct $\lambda_k\in \mathbb R$. Then $\,f$ is unbounded on $[0,\infty).$
\end{lemma}

\begin{proof} Let $T>0.$ Then 
 according to Corollary 2 of \cite{MV}, 
\[\int_0^T \big\vert    \sum_{k=1}^N c_k\, e^{i\lambda_k s}     \big\vert^2 ds \,=\, (T + 2\pi\theta\delta^{-1})\sum_{k=1}^N |c_k|^2\]
with a constant $\delta>0$ depending on the $\lambda_k$ and $|\theta|\leq 1$.
If $f$ were bounded on $[0,\infty), $ say $|f|\leq M$, this would imply that 
\[ \int_0^T \big\vert    \sum_{k=1}^n c_k\, e^{i\lambda_k s}     \big\vert^2 ds \, \leq M^2\int_0^T e^{-2as}ds\, \leq \frac{M^2}{2a},\]
a contradiction.
\end{proof}

Notice that only the first part of our proof of \ref{main_Theorem} uses uniform boundedness of hypergroup characters in order to settle boundedness of $F_\lambda$ in the case where $\, \text{Re}\,\lambda\,$ is contained in the boundary of $\text{co}(W.\rho).$
The rest  of the proof works  
equally for arbitrary root systems $R$ and arbitrary non-negative multiplicities $k\geq 0$, $k\not= 0,$ and the case $k=0$ is classical. Actually,  we have

\begin{corollary}
 Let $R\subset \frak a$ be an arbitrary root system, $k\geq 0$ a non-negative multiplicity function and $\rho = \rho(k).$ Then the associated
 hypergeometric function 
$t\mapsto F(\lambda,k;t)$ is unbounded on $\frak a$ if $\, \text{Re}\,\lambda\notin \text{co}(W.\rho).$ Moreover, 
$t\mapsto F(\lambda,k;t)$ is bounded on $\frak a$ if $\, \text{Re}\,\lambda\,$ is contained in the interior of $\text{co}(W.\rho).$
\end{corollary}

We return to our specific $BC$-cases and identify the dual space $\widehat{C_\mu}$ of the hypergroup $C_\mu$ with a subset of $\mathbb C^q$  via $\varphi_\lambda\mapsto \lambda$. Due to the condition $\,\overline \lambda \in\{w.\lambda, \, w\in W\}$ it is contained in the union of 
finitely many hyperplanes in $\mathbb C^q\cong \mathbb R^{2q}$ of real dimension $q$. Notice that the chamber $C$ is a proper subset of $\widehat{C_\mu}.$

The following is an immediate consequence of Opdam's Plancherel theorem (Thm. \ref{Plancherel}):

\begin{proposition} The Plancherel measure of  the hypergroup $C_\mu$ 
is given by the measure
\[ d\pi_\mu(\lambda)  = \frac{1}{|c(i\lambda, k_\mu)|^2}\,d\lambda.\]
on $\widehat{C_\mu}\subset \mathbb C^q$. Its support coincides with the chamber $C$. 
\end{proposition}

\section*{Acknowledgement} 

It is a pleasure to thank Angela Pasquale for fruitful discussions at an early stage of the work for this paper.

\end{document}